\pdfoutput=1
\documentclass[english]{IEEEtran}
\usepackage[T1]{fontenc}
\usepackage[latin9]{inputenc}
\usepackage{color}
\usepackage{mathtools}
\usepackage{amsmath}
\usepackage{bm}
\usepackage{amssymb}
\usepackage{graphicx}
\usepackage[final]{pdfpages}
\usepackage{cite}
\usepackage{mathtools,bbm}
\usepackage{babel}
\usepackage{algorithm}
\usepackage{algpseudocode}
\usepackage{pifont}
\usepackage{graphicx}
\usepackage{epstopdf}
\DeclareGraphicsExtensions{.eps}
\usepackage{verbatim}
\usepackage{esvect}
\usepackage{graphicx}
\usepackage{breqn}
\usepackage{subfig}
\usepackage{mathtools, cuted}
\usepackage{mathtools}
\usepackage{lipsum}

\newtheorem{definition}{Definition}

\begin{document}
\title{
Dynamic Formulation for Multistage Stochastic Unit Commitment Problem}
\author{Bita Analui and Anna Scaglione
\thanks{This research was supported in part by NSF Award No.1549923.}
\thanks{B. Analui and A. Scaglione are with the School of Electrical, Computer and Energy Engineering, Arizona State University, Tempe, AZ, USA. Emails: \texttt{\{bita.analui, anna.scaglione\}}@asu.edu}}
\maketitle

\begin{abstract} 
As net-load becomes less predictable there is a lot of pressure in changing decision models for power markets such that they account explicitly for future scenarios in making commitment decisions.  
This paper proposes to make commitment decisions using a dynamic multistage stochastic unit commitment formulation over a cohesive horizon that leverages a state-space model for the commitment variables. We study the problem of constructing scenario tree approximations for both original and residual stochastic process and evaluate our algorithms on scenario tree libraries derived from real net-load data.
\end{abstract}
\vspace{-4mm}
\section{Introduction and Motivation}\label{Motiv}
 \IEEEPARstart{T}{he} integration of 
 a greater percentage of renewable generation has prompted significant interest in incorporating 
stochastic optimization tools into the electric power market decision models. Perhaps, the most consequential of such decision models is the {\it Unit Commitment} (UC) problem \cite{Kazarlis&etal1996}, which determines which optimal set of generating units should be active to meet the electric net-load (load minus renewable generation) on the next hours or the next day. The net-load is unknown, but the UC schedule solves a deterministic problem using forecasts.  
Traditionally, uncertainty in power systems has been managed by scheduling for reserves to be online and then dispatching them in real time, to ensure that the system could sail through any single failure. This decision is referred to as the {\it Security Constrained Unit Commitment} (SCUC). Recently, there has been significant research focus on developing two stage stochastic SCUC models,  to account for uncertainty. The authors in \cite{Wu&etal2007} address the problem by modeling load uncertainties as well as the outages of generation units. In addition, in  \cite{Wang&etal2008} the SCUC model takes into account intermittent wind power. However, there is broad recognition that the SCUC decisions are insufficient to tackle the uncertainties from renewable generation, because contingencies are rare events, while random variations of the net-load will happen almost surely. If they are not accounted for, reserves can fall short, spillage of power becomes likely and a plethora of situations can render the real-time costs of operating the system very high, thwarting efforts to integrate more renewable power.
An alternative approach, introduced to mitigate the aforementioned inefficiencies, is the so called {\it Chance Constrained Unit Commitment} (CCUC), whose goal is to avoid conservative solutions that would lead to large power spillage \cite{Bouffard&etal2008} 
by allowing the violation constraints that correspond to random events based on a confidence level. An example is the work \cite{Wang&etal2014}, where the authors propose a UC formulation in which wind uncertainty is not modeled explicitly through a probabilistic model, but empirically through a set of actual past scenarios. The bundle of scenarios are used to produce a probability score for wind spillage and restrict it to the desired value through the decision variables.  
Whilst, the approach that we discuss in this paper, falls in a different class of decision models which capture the ensemble of behaviors that is expected directly, using a stochastic optimization framework \cite{Pflug&Pichler2014} and \cite{Ruszczynski2003}. Stochastic optimization in general replaces deterministic objectives with an expectations on the cost (or disutility) of the solution and in the constraints. It explicitly verifies that there is a feasible solution in all foreseeable scenarios. For UC problem,  this class of decision models are referred to as Stochastic Unit Commitment (SUC). To put our work into context, next we briefly review the most relevant literature on the SUC and then provide a synopsis of our main contributions relative to the state of the art. 
\textit{Stochastic Unit Commitment:} The first class (and the most commonly studied versions) of SUC formulations are two-stage stochastic optimization problems: 
the first stage decisions consist of choosing production units to meet the expected demand under operational costs and constraints (i.e. the commitment problem), while the second stage decisions consist of deciding the production level of the committed units to meet the actual demand under network constraints (i.e. the economic dispatch problem).
The model emulates the current decision model for short-term operation in power system that exhibits a two-stage structure. The wide range of existing literature in this category differ in algorithmic perspective, modeling and dimension of the uncertainty and incorporation of risk-averse objective function. The approach was originally proposed by \cite{Takriti&etal1996} and followed with wide range of scholarly works: \cite{Carpentier&etal1996} \cite{Lajda1981},\cite{Birge&Shiina2004}, \cite{Tuohy&etal2009}, \cite{Ruiz&etal2009}, \cite{Wang&etal2011} and the references therein.    
As opposed to the multistage stochastic unit commitment, two-stage approaches have been studied more significantly and more progress has been reported.  Even though multistage SUC (MSUC) approaches allow for smoother boundary conditions regarding the commitment variables (i.e. the \textit{on} and \textit{off} decisions) between consecutive decision horizons, it is often argued  in the literature that their complexity is prohibitive. 
It is true that scenario tree based multistage stochastic UC problems technically lead to large-scale linear mixed-integer problems which are hard to solve in general, but recent development have tried to curb such complexity. Specifically, the cutting plane algorithms proposed in \cite{Guan&etal2016} defines valid inequalities on scenario tree and leads to a more efficient and accurate solution compared to other approaches such as Progressive Hedging (PH)  \cite{Wets&etal2013} which as opposed to proven convergence results for convex problems, for non-convex cases, such as Mixed Integer Linear Programs, is only a heuristic algorithm. 

{\it Contributions}: In this paper we propose a novel dynamic MSUC (DMSUC) formulation which further improves the main benefit of the multistage approach through a state space formulation on cohesive instances of 24-hour horizon.  The additional advantage compared to the MSUC is that, as time goes on and new information becomes available, the decision maker is not stuck with commitment decisions that can still be reversed and that no longer reflect what is the best option at the current moment. The formulation requires to update the scenario tree dynamically, as well as to update the feasible action space due to the decisions that are made in the past and cannot be reversed. The benefit of our formulation is how the actions that are taken can be captured efficiently in {\it commitment state} variables. 

Noting that there are no performance guarantees in the construction of such 
scenarios in the context of Mixed Integer Stochastic Linear Programming (MISLP) solutions, we also discuss the issue of optimal quantization of actual scenarios bundles into the desired finite tree  approximation of the estimated filtration and explore the performance of the method, while constructing scenario trees based on actual net-load trajectories from PJM \cite{PJMload} and testing the DMSUC on the IEEE Reliability Test System (RTS) sample grid \cite{RTS1979}.

The remainder of this paper is organized in the following manner. In Section \ref{probStatement}, we first introduce the notations and mathematical concepts to set the layout for formulating DMSUC problem and discuss its performance in Sections \ref{DMSUC_Formula} and \ref{discuss_DMSUC}. In Section \ref{LoadScenario}, we discuss the issues arising when modeling the stochastic processes and present the algorithm to construct a library of scenario trees from net-load data. Our computational and ex-post performance  results are presented in Section \ref{ComRes}. Finally, Section \ref{Conclude} concludes the paper.       
\vspace{-4mm}
\section{Problem formulation}\label{probStatement}
\subsection{Dynamic Multistage Stochastic Optimization}\label{DMSO}
In the conventional Static Multistage Stochastic (SMS) optimization a decision maker minimizes the expected (constrained) cost of her actions under uncertainty, averaging over all possible realizations $\xi$ of a discrete time stochastic process $\Xi[k]$ over a single future horizon of $T$ time intervals that start at a certain time $k_0$. Without loss of generality, we can refer to the process $\Xi[k]$ for $k_0\leq k\leq k_0+T-1$ as: $\{\Xi[k_0+t]: t=0,\ldots,T-1\}.$
The scenario tree $\mathcal{T}=\{
\mathcal{V},\mathcal{E},\mathbb{P};\bm \xi\}$ represents the basic structure that is used in SMS optimization problems to represent the gradual unfolding of information through \textit{filtration} $\mathfrak{F}=\{\mathcal{F}_{k_0}\subset \mathcal{F}_{k_0+1}\subset \cdots\subset \mathcal{F}_{k_0+T-1}\}$ 
for the stochastic process $\{\Xi[k_0+t]: t=0,\ldots,T-1\}$. For now we shall assume that this mapping is bijective. However, later we will see in details that for computational purposes this mapping is reduced to be surjective, such that the scenario tree is seen as finite-discrete approximation of the true filtration. The tree is a directed graph, i.e. the paths go from the root to the leaves through edges $(u,v)\in \mathcal{E}$ where $u$ is connected to $v$, but not viceversa. 
Each node $v \in \mathcal{V}$ on the tree has a corresponding value $\xi_v\in {\bm \xi}$. The present, $\xi_0$, is unique (an event with probability \textit{one}) and represents the root of the tree. The collection of values of an individual path along the tree represent a possible realization of $\{\Xi[k_0+t]: t=0,\ldots,T-1\}$.  
Given any node $v$
The probability law ${\mathbb P}$ is specified by associating to each edge $(u,v)$ a probability $p_{u,v}$ as the conditional probability of the outcome $\Xi[k_0+t]=\xi_v$
given the unique realization ${\bm \Xi}_{t-1} = \bm \xi_{0:u}$ 
associated to the path that goes from the root to node $u$.
Using the chain rule the probability of a certain sample path or \textit{scenario} can be computed recursively as follows:
\vspace{-2mm}
\begin{align}
\pi_0&=Prob(\Xi[k_0]=\xi_0)=1\\
\pi_v&=Prob(\bm \Xi_t=\bm \xi_{0:v})=p_{u,v}\pi_u,~~(u,v)\in \mathcal{E}.
\end{align}
The general structure of the decision model becomes:
\begin{align}
\min_{\bm x}\mathbb E[C({\bm x},{\bm \xi})]~~
\mbox{s.t.}~~
{\bm x}\in {\mathcal X}({\bm \xi}),~~\bm \xi \sim \mathcal{T} \label{SMS}
\end{align} 
where both the convex cost function $C({\bm x},{\bm \xi})$ and constraints set ${\mathcal X}({\bm \xi})$ depend on the stochastic process outcomes.  The notation ${\bm \xi} \sim \mathcal{T}$ refers to the fact that the optimization provides solutions for all paths on the scenario tree. 
The components of decision array $\bm x$ themselves can be associated with filtration described by tree with exactly the same structure $\mathcal{T}_{opt}= \{\mathcal{V},\mathcal{E},\mathbb P ;\bm x\}$ but with the decision variables $x_v$ in lieu of the random outcomes $\xi_v$ associated to the nodes in $v\in \mathcal{V}$. This property is known as \textit{non-anticipativity} in stochastic optimization literature \cite{Rocka&wets1976}. The conventional formulation looks at a single period, but normally decisions of this kind have to be taken repeatedly in consecutive periods. This is certainly true in power market operations and  motivates our quest for a dynamic formulation. First, we separate the components of the decision array into the \textit{state} denoted by $\bm \sigma$, whose choice impact the future feasible decisions, 
and those that do not, denoted by ${\bm x}^{\prime}$, to partition vector $\bm x :=(\bm \sigma; {\bm x}^{\prime})$. Correspondingly, the variable associated to node $v\in\mathcal{V}$ is denoted by $ x_v =(\sigma_v;x^{\prime}_v)$. The state components $\bm \sigma$ are the part of decisions which constrain and couple consecutive decision horizons. Also note that, a dynamic formulation typically reduces this coupling to only two consecutive periods, hence the notion of {\it state}. 

Given all of the above, the dynamic formulation of an SMS  problem is conceptually a straightforward extension, since it simply amounts to changing the time index $k_0$ (which is entirely superfluous in the SMS) into the current time variable $k$.
Before, we formally define the Dynamic Multistage Stochastic (DMS) optimization problem, let us establish the following conventions:
\textbf{(i)} We will consider consecutive decision epochs of equal length $T$, each starting at present times $k$ for the corresponding horizon $\kappa$ with $k=\kappa T,~~ \kappa \in \mathbb{Z}$.     
\textbf{(ii)} $\mathcal{T}[\kappa]$ indicates the scenario tree that represents the 
filtration of the incoming segment $\Xi[\kappa T:(\kappa+1)T-1]$, conditioned on the past and present observations up to time $k$.\\
Mathematically, the dynamic nature of the DMS is captured by explicitly indicating that all structures and variables defined above depend on $\kappa$. Also, how the set of constraints in decision epoch $\kappa+1$ are affected by the actions and in particular the states of its prior epoch $\kappa$. 

Let us denote by $\mathcal{V}[\kappa]$ the set of all nodes on tree $\mathcal{T}[\kappa]$ and by $\mathcal{L}[\kappa] \subset \mathcal{V}[\kappa]$ the set of its leaf nodes. 
Let ${\bm x}^{opt}[k]$ be the optimum set of decisions that corresponds to the $\bm \xi \sim \mathcal{T}[\kappa]$. At time $(\kappa+1)T$ the decision maker will have observed the path that occurred in the previous period and, therefore, implemented a specific path of decisions out of
${\bm x}^{opt}[k]$, that correspond to this path.  Let this path correspond to the leaf that we denote as $v^*\in \mathcal{L}[\kappa]$, the decisions that were taken be ${\bm x}^{opt}_{0:v^*}[k]$ and, in particular, the state at the leaf be $\sigma_{v^*}^{opt}[k]$.
Only this state $\sigma_{v^*}^{opt}[\kappa]$ carries the necessary information about the past decisions that need to be transferred in to the next  horizon. 
Therefore, the optimization problem for epoch $\kappa+1$ is an extension of \eqref{SMS} and formulated as: 
\begin{align}
\min_{\bm x[\kappa+1]} &~~\mathbb E[C_{\kappa+1}({\bm x}[\kappa+1],{\bm \xi}[\kappa+1])]\label{DMS_obj}\\ 
\mbox{s.t.}&~~{\bm x}[\kappa+1]\in {\mathcal X}({\bm \xi}[\kappa+1],\sigma_{v^*}^{opt}[\kappa]),
\label{DMS_feasibleSet}\\
 &~~\bm \xi[\kappa+1] \sim \mathcal{T}[\kappa+1] \label{DMS}
\end{align}
Further in Sections \ref{TreeLib} and \ref{ComRes} we will explain in greater details how $\mathcal{T}[\kappa+1]$ is chosen and the state of the system from $\mathcal{T}_{opt}[\kappa]$ is transferred into it.
\vspace{-4mm}
\subsection{Dynamic Multistage Stochastic Unit Commitment}\label{DMSUC_Formula}
In this section we first introduce a new nodal formulation for MSUC that makes the generalization to our Dynamic MSUC (DMSUC) decision algorithm straightforward. Prior to that, all the component of a MSUC problem and its nomenclature is defined for a static formulation. 
\subsubsection{Static formulation of the MSUC}
In the MSUC we consider the underlying scenario tree $\mathcal{T}$ representing the probabilistic information about the net-load. For now we consider the parameters that define the cost function as deterministic, although it is straightforward to extend what follows to the case where also the future cost parameters are unknown and, therefore, describe them in the same fashion as stochastic parameters on the scenario tree. \\
The power system is composed of several generating units $g$ at each bus $ b\in \mathcal{B}$. All the buses are connected through a set of transmission lines $(l,l^{\prime})\in \mathcal{B}\times\mathcal{B}$. The set of generating units at bus $b$ is denoted by $\mathcal{G}^b$ and their generation schedule corresponding to a certain node $v$ on the scenario tree is the pair of  continuous and binary variables $(x^g_{v},y^g_{v}) \in \mathbb{R}_+  \times \mathbb{B} $ with $g\in \mathcal{G}^b$ and $v\in \mathcal{V}$. In addition to the typical schedule and commitment variables $x_v^g \in \mathbb{R}_{+}\  \text{and}\  y_v^g\in \mathbb{B}$, two binary variables $\overline{s}^g_v$ and $\underline{s}^g_v$ indicate the switching action from \textit{off} to \text{on} and from \textit{on} to \text{off} respectively. 
The switching variables 
will take respectively value 0 unless the unit is turned on or off, in which case they will take respectively values 1. However, it is shown in subsequent sections that variables $\overline{s}^g_v$ and $\underline{s}^g_v$ are expressed in terms of binary variable $y^g_v$. Therefore, they can be relaxed to belong to $[0,1]$.\\ 
In addition to these standard variables, we introduce two {\it state} variables to handle:
\begin{itemize}[noitemsep,topsep=0pt]
\item minimum-up: $O_n^g$ 
\item minimum-down: $O_f^g$ 
\end{itemize}
time of a unit $g \in \mathcal{G}^b$, called $o^g_v$ and $d^g_v$.\\ 
The variable $o^g_v,\ v \in \mathcal{V}$, for $v$ corresponding to $t^{th}$ time instant in the horizon, is the residual time unit $g$ needs to stay on after $t$, which depends on the parent state $o^g_{v_-}$; so, only when $o^g_v=0$ the unit can be turned off and the state persists for the next generations as long as the machine continues to stay on or, if is switched off, for as long as it is off and not switched on again. The variable $d^g_v, v\in {\mathcal V}$, is the complementary variable handling the off time constraint, and represents the residual time the unit needs to remain off after time $t$ depending on the state of the parent of node $v_{-}$. Similarly, only when $d^g_v=0$ the optimization can turn the unit on and $d^g_v=0$ persists after the unit is turn on until the next switch-off event.
Note that this set of constraints are applied only to those units $g\in \mathcal{G}^{s}\subset \bigcup_{{b\in\mathcal{B}}}\mathcal{G}^b$ with $O_n^g > 1$ and $O_f^g > 1$.
Below the full nodal formulation of MSUC  problem is presented and in the following the corresponding constraints are described.
\vspace{-1.5mm}
\begin{subequations}
\begin{align}
\min & ~\sum_{v\in {\mathcal V}}\pi_v \sum_{g\in \mathcal{G}} [S^g_v+C^g_v+\frac{o_v^g}{O_n^g}+\frac{d_v^g}{O_f^g}] \label{cost}\\ 
\text{w.r.t} &~~~(\mathbf{y},\mathbf{o}, \mathbf{d}, \mathbf{x}, \mathbf{\overline{s}}, \mathbf{\underline{s}}) \notag\\
\text{s.t.} &~~~ {\bm \xi} \sim \mathcal{T} \label{N.A.}
\\&~ \sum_{b\in\mathcal{B}}\bigg(
\sum_{g\in \mathcal{G}^b}x^g_v-\xi^b_v\bigg)=0 \label{balance} ~~~~~~~~~~~~~~~~~ (\forall v\in \mathcal{V})\\
&~-L_{ll^{\prime}}\leq 
\sum_{(l,l^{\prime})\in\mathcal{B}\times\mathcal{B}}D^b_{ll^{\prime}}\bigg(
\sum_{g\in \mathcal{G}^b}x^g_v-\xi^b_v
\bigg)\leq L_{ll^{\prime}} \label{lineFlow}\\
&~ o^g_v \geq \overline{s}^g_v(O_n^g-1)\label{minUp1}\\
&~ \max\{0,o^g_{v_-} -y^g_v\}\leq o^g_v\leq o^g_{v_-}+\overline{s}^g_v (O_n^g-1)\label{minUp2}\\
&~ o^g_{v_-}-o^g_v \leq y^g_v\leq 1\label{minUp3}~~~~~~(\forall v\in \mathcal{V}\smallsetminus\{0\},~g\in\mathcal{G}^{s})\\
&~ d^g_v \geq \underline{s}^g_v(O^g_f-1)\label{minDown1}\\
&~\max\{0,d^g_{v_-}-1+y^g_v\}\leq d^g_v\leq d^g_{v_-}+\underline{s}^g_v (O^g_f-1)\label{minDown2}\\
&~0\leq y^g_v\leq 1-d^g_v+d^g_{v_-}\label{minDown3}~(\forall v\in \mathcal{V}\smallsetminus\{0\},~g\in\mathcal{G}^{s})\\
&~y^g_v-y^g_{v_-}\leq \overline{s}^g_v \label{on}\\
&~\underline{s}^g_v = y^g_{v_-}- y^g_v + \overline{s}^g_v \label{off}~~~~~~~~~~~~~~~~(\forall v\in \mathcal{V}\smallsetminus \{0\})\\
&~\underline{G}^gy^g_v\leq x^g_v\leq \overline{G}^gy^g_v \label{production}~~~~~~~~~~~~~~~~~~~~~~~(\forall v\in \mathcal{V})\\
&-{{\underline{G}^g}^{\prime}}y^g_v\leq x^g_v-x^g_{v_-}\leq {{\overline{G}^g}^{\prime}}y^g_v  \label{ramping}~~~(\forall v\in \mathcal{V}\smallsetminus\{0\})\\
&(y^g_{0-},o^g_{0-},d^g_{0-})=(0,0,0) \label{state}\\
& \mathbf{y}\in \mathbb{B}^{|\mathcal{G}|\times|\mathcal{V}|},\ \mathbf{o}, \mathbf{d}, \mathbf{x} \in \mathbb{R}_{+}^{|\mathcal{G}|\times|\mathcal{V}|},\  \mathbf{\overline{s}}, \mathbf{\underline{s}}\in [0,1]^{|\mathcal{G}|\times|\mathcal{V}|} \notag 
\end{align} 
\end{subequations} 
In the formulation above, \eqref{cost} expresses objective i.e., to minimize the expected total cost. The actions of switching on and off the units come with the so called start up $\overline{S}^{g}$ and shut down cost $\underline{S}^{g}$. In addition to these costs, there is a price bid that is a piece-wise linear function of the production level $x^g_{v}$. 
Therefore the nodal cost term is expressed as:
\vspace{-1.75mm}
\begin{align}
S^g_v= \overline{S}^{g}\overline{s}^g_v+\underline{S}^{g}\underline{s}^g_v \\
C^g_v= c_1^g x^g_{v}+ c_0^g y^g_{v}.
\end{align}
\eqref{balance} and \eqref{lineFlow} constraints are added to account for the grid balance requirement under load uncertainty  and line flow constraints for flow limits $[-L_{ll^{\prime}},L_{ll^{\prime}}]$ for each transmission line $(l,l^{\prime})\in \mathcal{B}\times \mathcal{B}$ connecting bus $l$ to bus $l^{\prime}$. 
Constraints \eqref{minUp1}-\eqref{minUp3} will enforce the behavior which was described in defining the $o^g_v$ for all $v\in \mathcal{V}$ which have as parent node (immediate predecessor) $v_{-}$. Constraint \eqref{minUp1} sets the state variable to the residual on-time when the unit is switched on.
Constraint \eqref{minUp2} together with $\frac{o^g_v}{O^g_n}$ term in the objective, guarantee that the unit has a monotonically decreasing state when $1\leq o^g_{v_-}\leq O^g_n-1$. But once the unit reaches $o^g_{v_-}=0$ it keeps the unit in that state until $\overline{s}^g_{v_{+}}=1$ (by convention, $v_{+}$ denotes an arbitrary successor of node $v$ or $v\prec v_{+}$). Note that, in this case the state goes from 0 to $O^g_n-1$ so is not monotonically decreasing and that is why the upper-bound for $o^g_v$ contains the term $\overline{s}^g_v(O^g_n-1)$ which allows the transition to a higher state to happen.
Constraint \eqref{minUp3} forces the unit to remain on (i.e. $y^g_v=1$) as long as the state is monotonically decreasing by one, except when the parent $o^g_{v_-}=0$, which means that the unit cannot be turned off until the minimum-up time is over. Note also that, when the unit is off $y^g_v=0$, irrespective of what is going on with $d^g_v$, the state must be $o^g_v=0$ and remains unchanged until the unit is switched on. 
The explanation of the constraints in  \eqref{minDown1}-\eqref{minDown3} follow exactly the same line of reasoning. Constraints \eqref{on}-\eqref{off} express the switching variables in terms of the commitment variable. Moreover, \eqref{production}-\eqref{ramping} guarantee the production and ramping capacity limits respectively. 
It is also essential to address that the non-anticipativity constraint is already inherent in the nodal formulation.   
\subsubsection{Dynamic formulation of the MSUC}\label{DMSUC}
After one update in the tree 
${\mathcal{T}}[\kappa]$, as we discussed for the general case in Section \ref{DMSO},
going from the static to the dynamic formulation amounts to finding how the past (optimal) decision variables constrain the decisions in the decision horizon corresponding to scenario tree $\mathcal{T}[\kappa+1]$. 
The partition for the decision variables of MSUC problem into state variables and the rest is as follows:
\vspace{-2mm}
\begin{align*}
\underbrace{(y^g_{v},o^g_{v}, d^g_{v}}_{\sigma_v};\underbrace{x^g_{v}, \overline{s}^g_{v}, \underline{s}^g_{v})}_{x^{\prime}_v}.
\end{align*}
Recall \eqref{DMS_obj}-\eqref{DMS} and as we stated in the previous section, let $v^*\in \mathcal{L}[\kappa]$ be the leaf node indicating the path that was observed during the $\kappa^{th}$ horizon and just elapsed. In addition, 
let $(y^g_{v^{*}}[\kappa],o^g_{v^{*}}[\kappa],d^g_{v^{*}}[\kappa])^{opt}$ be the triplet of optimum decisions that correspond to the leaf node $v^*$. To compute the optimal decisions corresponding to $\kappa+1$ epoch, with the dynamic formulation, the state value in \eqref{state} must be initialized with:
 \begin{align}
(y^g_{0-},o^g_{0-},d^g_{0-})= (y^g_{v^{*}}[\kappa],o^g_{v^{*}}[\kappa],d^g_{v^{*}}[\kappa])^{opt}.
 \end{align}
After the state of the system is transferred from decision horizon $\kappa$ to $\kappa +1$, the MSUC problem \eqref{cost}-\eqref{state} is solved for underlying tree $\mathcal{T}[\kappa+1]$. By finding solutions for all nodal variables
$(\mathbf{y}[\kappa+1],\mathbf{o}[\kappa+1], \mathbf{d}[\kappa+1], \mathbf{x}[\kappa+1], {\overline{\mathbf s}[\kappa+1]}, {\underline{\mathbf s}[\kappa+1]})$ the decision maker performs the action 
and awaits the information regarding the next period and the process continues for every incoming interval.
\vspace{-4mm}
\subsection{Discussion on the DMSUC complexity and performance}\label{discuss_DMSUC}
The benefit of the DMSUC is that, it allows, to an extent, to co-optimize the commitment schedule and real time balancing. It also enables the market players to change their positions dynamically, increasing the economic competition. Also, while it is true that one needs to solve many optimization problems, it is possible that a relatively small look-ahead is sufficient to provide schedules that are more efficient and stable compared to the state of the art. Meanwhile, ideally for a scenario tree to guarantee the optimality of the decisions with respect to the true model, it should encompass the totality of all possible outcomes of the underlying stochastic process. To this end, unless the distribution of stochastic process has finite and countable outcomes, no scenario tree with finite dimensions can completely represent the stochastic process and its corresponding probability model. In addition, even in presence of finite and countable probability space, we know that with increasing number of decision stages the number of scenarios grow exponentially. Hence, solving a day ahead stochastic unit commitment problem in a in a $24$-hour horizon can be daunting, but our proposed approach exploits the quality of decision along the whole horizon ($24$-hours) while solving a series of the DMSUC problem over shorter periods and transferring the optimal state of the system in to the next period with the arrival of the new observation.
In the next section we complement our dynamic formulation with a dynamic construction of an approximate scenario tree that is also dynamic and data driven.
\vspace{-1mm}
\section{Dynamic Construction of Load Scenario Trees}\label{LoadScenario}
The most popular and widespread numerical solution method for multistage stochastic optimization is employing a finite and discrete approximation of the true probability model of stochastic process in the form of scenario tree. Before 
we introduce our method, we discuss next the type of errors expected in our design. 
\vspace{-4mm}
\subsection{Modeling and Approximation Error}
In formulating the MSUC and its extension to DMSUC, we used notation $\mathcal{T}[\cdot]$ to introduce the \textit{true} scenario tree which includes the true probability model and filtration. In practice, however, capturing the true underlying model is unlikely since:\\ 
{\bf (i)} In most cases the true probability model of the stochastic process is not available. Therefore,
one can rely on past observations to fit a parametric model or infer the empirical model in a non-parametric fashion. In both cases \textit{modeling} errors are consequential.\\
{\bf (ii)} Even if the true probability model is available, when the process is continuous, for computational purposes, a finite approximation of $\mathcal{T}[\cdot]$, denoted by $\widetilde{\mathcal{T}}[\cdot]$, is required. 
The approximation, to the extent that is computationally tractable, will cause the \textit{approximation} error in the scenario tree of actions (or decisions).\\
These two facts result in errors that are compounded and therefore, sub-optimal solutions. Ignoring the modeling errors, it is still unknown how to minimize the deviation form the optimal objective of the real problem under a constrained structure for the approximate tree.  However, it is shown in \cite{Pflug&Pichler2012} that under specific continuity assumptions, the approximation error which is measured directly by the distance between the true process and the approximate one can provide upper bound information about the performance of the optimization. Unfortunately, with the inclusion of binary variables in the MSUC (and DMSUC) such results are not applicable. To the best of our knowledge, the issue of ensuring that approximations errors on the scenario tree result in minimal errors in the integer variables has not been specifically addressed in the literature. Only in \cite{Maggioni&etal2016} the authors propose general bounding techniques based on solving independent sub-problems on scenario tree that is valid for mixed-integer MSO problems too.
What can be certainly said is that the bounds proposed in \cite{Pflug&Pichler2012} would still be applicable if and only if the \textit{optimal} integer variables remain unchanged in the surjective mapping between each 
path in $\mathcal{T}[\cdot]$ and the corresponding path in $\widetilde{\mathcal{T}}[\cdot]$. 
In fact, if the original problem on 
$\mathcal{T}[\cdot]$ corresponds to a set of optimal integer variables $\mathbf{y}^{opt}$, 
for a certain set of children at stage $t$ the optimum decision tree 
includes $m_{t}$ distinct choices for the integer variables
$y_{v_t}^{opt}$, then 
$\widetilde{\mathcal{T}}[\cdot]$ must have at least $m_{t}$ children at that stage for each node $v_t$. In general $m_{t}$ will be quite a bit smaller than the exponentially large possibilities that exist for each 
$y_{v_t}^{opt}$, also because of the constraints. 
We leave the theoretical study of bounds for the error in mixed-integer formulations as future work and use heuristics to construct our approximation whose efficacy is only numerically tested.       
\vspace{-3mm}
\subsection{Scenario Tree Construction}
As mentioned in Section \ref{Motiv}, we aim to model the uncertainty through a library that contains load scenario trees. These trees can be used as the underlying structures to DMSUC formulation. As explained in Section \ref{DMSUC} each tree is selected to represent the uncertainty in the incoming horizon.
One important aspect seldom mentioned is that, as long as a reliable forecast method is available that results in a residual error with  smaller variance than the original process, this residual error process of the forecast fully represents the future uncertainty. The forecast of the value $\Xi[k]$ is a function $h(\cdot)$ of past samples of the process. Looking at all possible realizations for this set of past samples, the estimate $\hat{\Xi}[k]$ is itself a random process:
\begin{equation}
\hat{\Xi}[k]=h(\Xi[k-1],\Xi[k-2],\ldots).
\end{equation}
Rather than trying to obtain an approximate tree for the process $\Xi[k]$, given an admissible estimator $\hat{\Xi}[k]$ of $\Xi[k]$, a more sensible option is to create scenario trees for the residual error process
$Z[k]=\Xi[k]-\hat{\Xi}[k]$. Interestingly, at time $k$, we can view the process as the result of the following dynamics:
\vspace{-2mm}
\begin{equation}
\Xi[k]=h(\xi[k-1],\xi[k-2],\ldots)+Z[k],
\end{equation}
and view the sequence used for the forecast $\beta[k-1]=(\xi[k-1],\xi[k-2],\ldots)$ as the state of the dynamics at time $k$, where the state transition between two consecutive states is:  
\vspace{-2mm}
\begin{equation}
\beta[k-1]\mapsto\beta[k]=(\xi[k],\beta^{-}[k-1]),
\end{equation}
with $\beta^{-}[k-1]$ be vector $\beta[k-1]$ deprived of its last entry. 
Having formed a filtration for the residual process:
\begin{equation}
{\bm \zeta}[\kappa]\sim \mathcal{T}_Z[\kappa] \label{RezTree}
\end{equation}
to capture the uncertainty in the dynamics over the given horizon, we can obtain the corresponding nodal values for $\xi_v$ as follows:
\vspace{-2mm}
\begin{align}
\xi_v[\kappa]=h(\beta_{v_-}[\kappa])+\zeta_v[\kappa],~~
\beta_{v}[\kappa]=(\xi_v[\kappa],\beta^{-}_{v_-}[\kappa])\label{DynamicEq}
\end{align}
and the state $\beta_0[\kappa]=(\xi_0[\kappa],\beta^{-}[\kappa T-1])$, 
where $\beta[\kappa T-1]$ is the state that was actually observed at time $\kappa T-1$ just before the start of the current decision horizon. Thus, expression \eqref{RezTree} replaces constraint \eqref{N.A.} in optimization problem \eqref{cost}-\eqref{state} and $\xi_{v},  \forall v\in \mathcal{V}$ is explicitly derived from \eqref{DynamicEq}.
This approach may easily be utilize to model the finite difference process of the net-load, which expresses the ramping process, to significantly improve the approximation error.    
\vspace{-2mm}
\subsection{Exploiting Cyclostationarity Property of Load Data}
\begin{figure}
\centering
\includegraphics[width=0.49\textwidth,height=4cm]{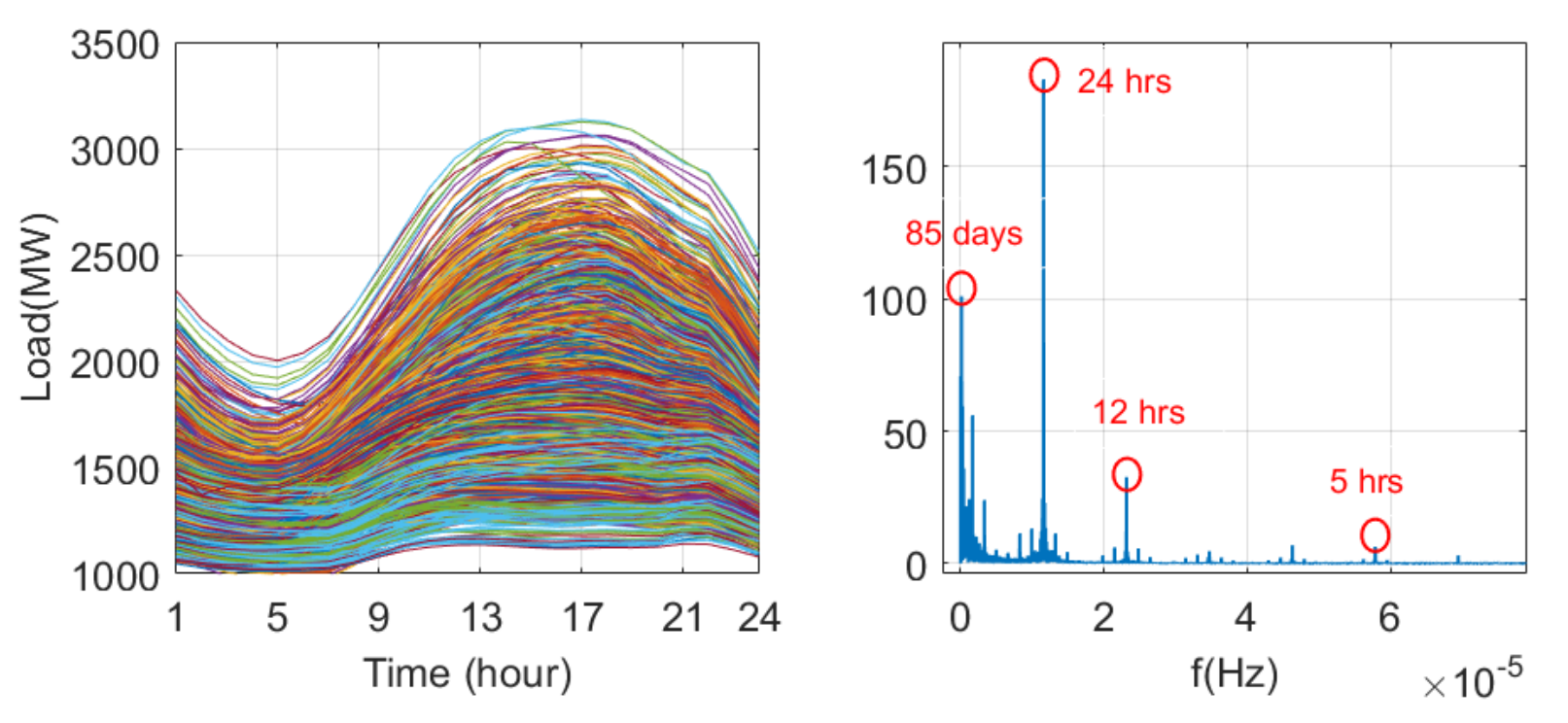}
\caption{\footnotesize{Left: Hourly Summer Load. Right: Cyclic Frequencies}}\label{load_fft}
\end{figure}

\begin{figure}
\begin{minipage}[t]{.48\textwidth}
\centering
\subfloat[Load Trajectories (MW) correspond to quantization bins]
{\includegraphics[width=.98\linewidth,height=5cm]{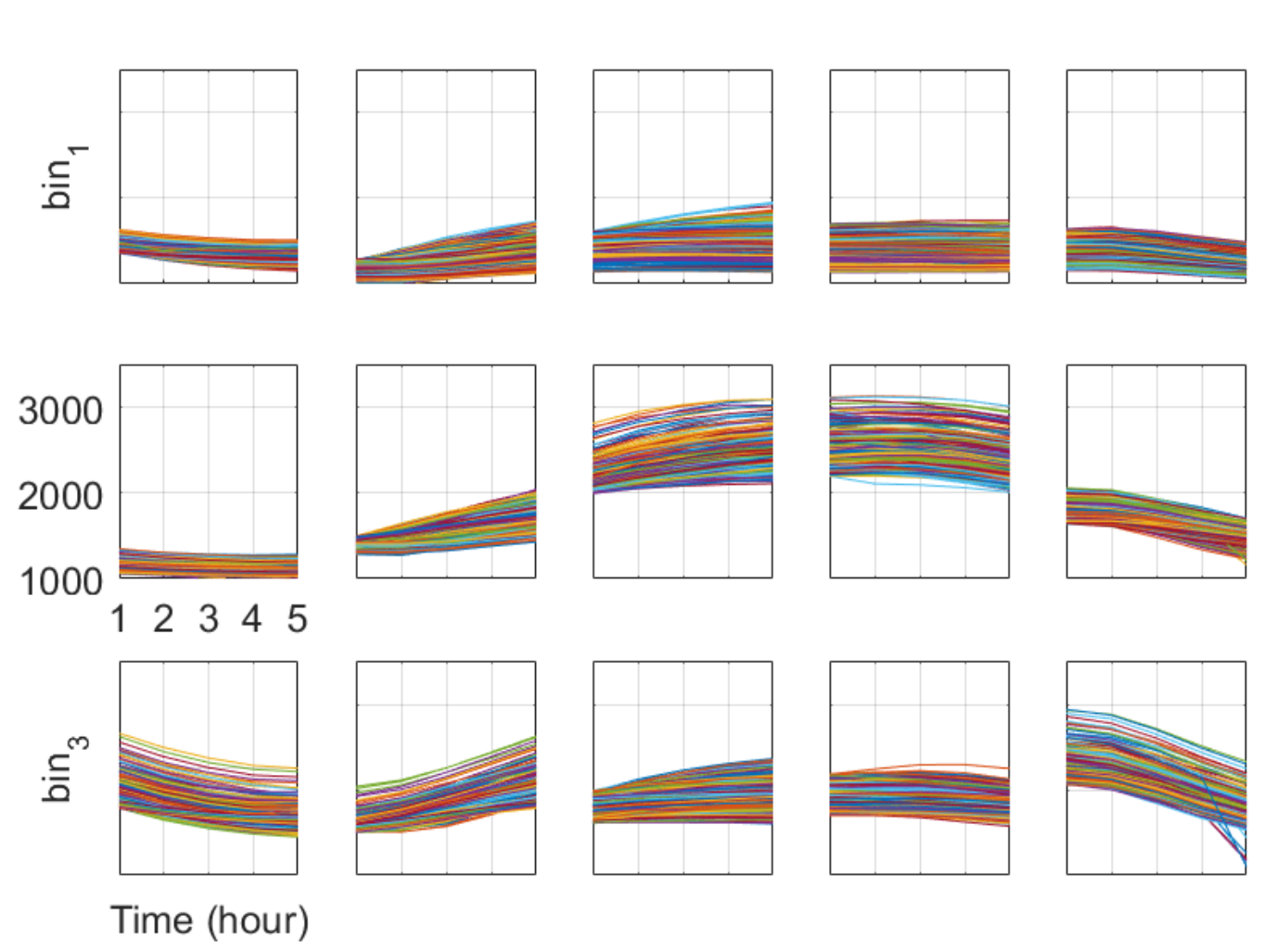}}\\
\subfloat[Scenario Trees corresponding to Load quantization bins]
{\includegraphics[width=.98\linewidth,height=5cm]{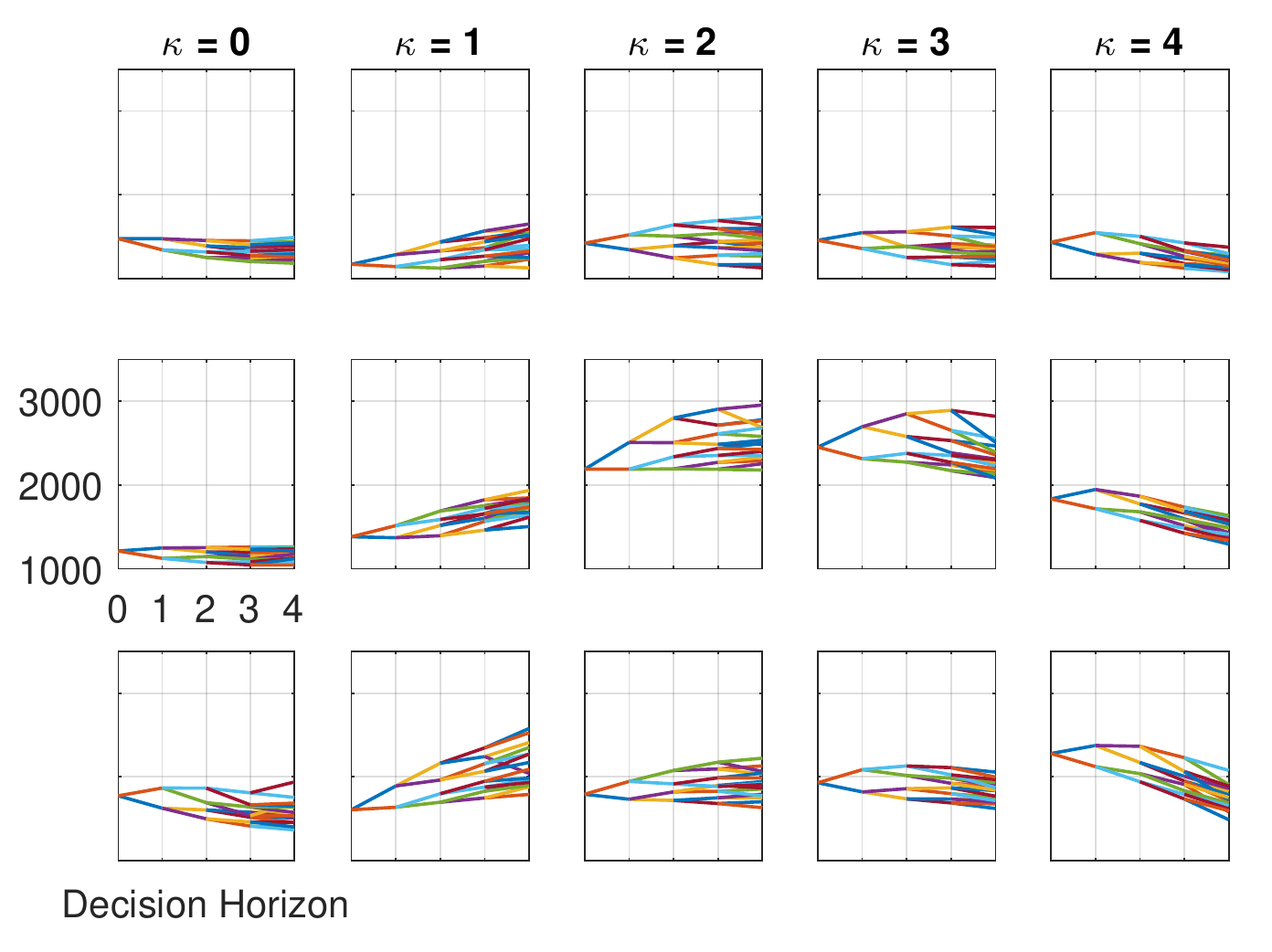}}\\
\caption{\footnotesize{(a) Load (MW) trajectories  
$\xi^l_{\kappa}$. (b) Load scenario trees from 00:00 to 01:00AM+1, with $c=3$ quantization bins and $T=5$. (For all figures ylim = [1000, 3500] Megawatt.)}}\label{BucketTree}
\end{minipage}\hfill
\end{figure}

To illustrate the proposed procedure we use hourly measurements of PJM \cite{PJMload} load data, during summer months for years 2000-2015 from the different load zones within the PJM regional transmission organization. 
As an illustration, the original summer load trajectories are depicted in Fig.\ref{load_fft}, for years 2000-2015. Obviously, the load trajectories show cyclic  features. It is evident from Fig. \ref{load_fft}, that load is not a stationary process and exhibits interdaily, daily and annual seasonalities.
Therefore, hourly summer load can be treated as a cyclostationary stochastic process with period equal to one day (24-hours) and consecutive periods of 5-hours within the day. More formally, a cyclostationary process is defined.
\begin{definition}
Discrete time stochastic process $\Xi[k]$ is cyclostationary in the strict sense with period $P$ if and only if for any set of indexes $(k_1,\ldots,k_n)$
the vector $(\Xi[k_1],\ldots,\Xi[k_n])$ and $(\Xi[k_1+\ell P],\ldots,\Xi[k_n+\ell P])$ have the same joint distribution for any $\ell\in \mathbb{Z}$.
\end{definition}
In the subsequent part of this section, we present the procedure for constructing a library of net-load scenario trees from empirical samples.
By employing the assumption that hourly summer load data exhibits cyclostationary property with statistics repeating every $P$ period, the sample paths from which the trees are constructed, are taken from records of the same cycle for the corresponding decision horizon $T$.\\
Let us denote the discrete time stochastic load process by $\Xi$ with $\xi$ be its realizations and the decision horizon for which the scenario tree is constructed by $T$. 
In addition, $\kappa \in \mathbb{Z}$ indexes the arbitrary consecutive decision epochs of length $T$. As it was addressed in \cite{Analui&Scaglione2016}, it is essential to account for all possible scenario trees that carry distinct information about the present. Therefore, prior to scenario tree generation, we first perform \textit{root} node quantization procedure by employing Lloyd-Max clustering algorithm (cf.\cite{Lloyd1982}). The optimization procedure is done in order to find positions or centroids  $\tilde{\xi}^l_0[\kappa T]$ for $l=1,\ldots, c$ number of clusters and decision epoch $\kappa$. 
Each of these optimal solutions ${\tilde{\xi}}_0^{l}[\kappa T], l=1,\ldots, c$ is a centroid associated with a \textit{Voronoi} cell $V^l_{\kappa}$ which contains particular and not necessarily equal number of load trajectories. 
Having determined the root nodes and corresponding sample trajectories as shown in Fig.\ref{BucketTree}-(a), next we will address the scenario tree generation through stochastic approximation approach.
\vspace{-5mm}
\subsection{Library of Scenario Trees}\label{TreeLib}
A stochastic approximation algorithm as proposed in \cite{Pflug&Pichler2015} is used in order to construct the library of scenario tree: $\mathfrak{T}^{0} = \{ \widetilde{\mathcal{T}}^l[\kappa]:\kappa \in \mathbb{Z} \  \text{with root}\ \  {\tilde{\xi}}_0^{l}[\kappa T], l=1, \ldots, c\}.$ 
Fig.\ref{BucketTree}-(b), depicts the library $\mathfrak{T}^{0}$ for $c=3$. Superscript $0$ corresponds to the fact that library $\mathfrak{T}$ is constructed base on the historical data, however the library can be dynamically updated and evolve with the arrival of new information. 
\\Let us fix the graph of the tree $\{\mathcal{V}, \mathcal{E}\}$ to a \textit{binary} structure (i.e., $\forall \  v\in \mathcal{V}, |v_{+}|=2$) in order for the filtration to be represented uniformly along all consecutive decision horizons. The general part of the algorithm, which is the stochastic approximation iteration, is explained below and the construction of the library of scenario trees is briefed in Algorithm \ref{Library_of_Trees_Alg.}.\\ 
Let $\tilde{{\bm \xi}}^{(0)}$ be the initialization of the nodes based on the empirical marginal distributions (cf. \cite{Pflug&Pichler2015}, Algorithm 3-(i)). Trajectory:\vspace{-2mm} \[tr^{(n)}(\xi^l)=[\tilde{\xi}^l_0[\kappa T], \xi^l[\kappa T-\ell P+1:(\kappa+1)T-\ell P-1]],\] with $n=1,\ldots, N$ denoting the number of iterations, is randomly picked from the sample set within each quantization Bin. 
By searching stage-wise through the paths on the tree, the candidate path for stochastic approximation step is found. We introduce set $\mathcal{P}_0\times\mathcal{P}_1\times\mathcal{P}_2\times\ldots\times\mathcal{P}_{T-1}$ and it elements by $\rho$. Each vector $\rho$ contains the node indices on the tree corresponding to each scenario path from root node to the leaves. In addition, a counting variable $\iota_{\rho}^{(0)}= \vv{0}$ is initialized to $zero$ and assigned to each scenario path $\rho$. 
We find $\rho^{*} \in \mathcal{P}_0\times\mathcal{P}_1\times\mathcal{P}_2\times\ldots\times\mathcal{P}_{T-1} $ such that: 
 \begin{equation}
 \tilde{{\bm \xi}}_{\rho^{*}}^{(n-1)} \in \mbox{argmin} ~\mathbf{d}(tr^{(n)}(\xi^{l}), \tilde{{\bm \xi}}_{\rho}^{(n-1)}), \label{disT}
 \end{equation}
 $\mathbf{d}(\cdot,\cdot)$ is an $l_2$-norm here. Consequently, path $\tilde{{\bm \xi}}_{\rho^{*}}^{(n-1)}$ is updated through the following gradient  equation:
\begin{align}
\tilde{{\bm \xi}}_{\rho^{*}}^{(n)} = &~\tilde{{\bm \xi}}_{\rho^{*}}^{(n-1)}-2a_{n} \mathbf{d}(tr^{(n)}(\xi^{l}), \tilde{{\bm \xi}}_{\rho^{*}}^{(n-1)})\triangledown \mathbf{d}(tr^{(n)}(\xi^{l}), \tilde{{\bm \xi}}_{\rho^{*}}^{(n-1)}) \label{update}\\
\vspace{-5mm}
\iota^{(n)}_{\rho^{*}} =&~ \iota^{(n-1)}_{\rho^{*}}+\mathbbm{1}, 
\end{align}
with $a_n$ being the step size (cf. \cite{Pflug&Pichler2015}, Remark. 10). Iterative steps \eqref{disT} and \eqref{update} are repeated for $n=1,\ldots,N$ and until all the successor nodes on the tree are fixed. Finally, the conditional probabilities $p_{\rho}$ are computed:
\vspace{-2mm}
\begin{equation}
 p_{\rho}=\frac{\iota^{(N)}_{\rho}}{N}, ~~ \forall \rho \in~~ \mathcal{P}_0\times\mathcal{P}_1\times\mathcal{P}_2\times\ldots\times\mathcal{P}_{T-1}  \label{cProb}
\end{equation}
The following algorithm summarizes the procedure for generating the scenario tree library $\mathfrak{T}^0$.
\begin{algorithm}
  \caption{Constructing Library of Scenario Trees from Empirical Load Samples}\label{Library_of_Trees_Alg.}
\footnotesize{\begin{algorithmic}[1]
    \State $\{\mathcal{V}, \mathcal{E}\}$, $\tilde{{\bm \xi}}^{(0)}$,   $c$, $\kappa$, $N$, and $\mathfrak{T}^0=\{\}$ are given. 
    \State $\forall  \kappa$ and for $l=1,\ldots,c$, centeroids ${\tilde{\xi}}_0^{l}[\kappa T]$ are computed.
    \State repeat for $\kappa$  \label{Horizon}
    \State repeat for $l$ \label{Bin}
    \State for $n=1:N$  \label{Iter}
    \State randomly pick trajectory $tr^{(n)}(\xi^{l})$ from corresponding bin
    \State find $\rho^{*}$ such that: 
    $\tilde{{\bm \xi}}_{\rho^{*}}^{(n-1)} \in \mbox{argmin} ~\mathbf{d}(tr^{(n)}(\xi^{l}), \tilde{{\bm \xi}}_{\rho}^{(n-1)})$
    \State update the path:
 \[\tilde{{\bm \xi}}_{\rho^{*}}^{(n)}=\tilde{{\bm \xi}}_{\rho^{*}}^{(n-1)}-2a_{n} \mathbf{d}(tr^{(n)}(\xi^{l}), \tilde{{\bm \xi}}_{\rho^{*}}^{(n-1)})\triangledown \mathbf{d}(tr^{(n)}(\xi^{l}), \tilde{{\bm \xi}}_{\rho^{*}}^{(n-1)})\]
    \State update $\iota^{(n)}_{\rho^{*}}$
   \State $n = n+1$ and return to \eqref{Iter} 
    \State $\widetilde{\mathcal{T}}^l[\kappa] = \tilde{\bm \xi}^{(N)}$, add new tree to the library: $\mathfrak{T}^0=\mathfrak{T}^0 \cup \widetilde{\mathcal{T}}^l[\kappa]$  
  \end{algorithmic}}
\end{algorithm}\\
Since tree library $\mathfrak{T}^0$ is constructed based on the historical load data it can be viewed as a repository of trees. 
Obviously, considering that scenario tree construction is purely data driven, it is substantial to update the trees in the library with observation of a complete path. This update will relocate the centroids in order to account for new realized path.   
\subsubsection{Here-and-now update}
The root node on the scenario tree contains the information about the present event with probability \textit{one}. Therefore, when the $\mathrm{p}$resent value $\tilde{\xi}_{\mathrm{p}}[\kappa T]$ of the load for decision epoch $\kappa$ is known, it is used to map library $\mathfrak{T}^0$ into the correct $\widetilde{\mathcal{T}}^l[\kappa]$ by finding index $l^{*}$ such that:
 \begin{equation}
{\tilde{\xi}}_0^{l^*}[\kappa T] \in \mbox{argmin}_l ~\mathbf{d}(\tilde{\xi}_{\mathrm{p}}[\kappa T],{\tilde{\xi}}_0^l[\kappa T])
\label{pickingTrees}
\end{equation}
When index $l^*$ is found, tree $\widetilde{\mathcal{T}}^{l^*}[\kappa]$ is picked from the library and its root node is replaced with $\mathrm{p}$resent value $\tilde{\xi}_{\mathrm{p}}[\kappa T]$. The corresponding tree is also denoted by the augmented tree
$\widetilde{\mathcal{T}}^{l^*}_{a}[\kappa] := (\widetilde{\mathcal{T}}^{l^*}[\kappa];\tilde{\xi}_{\mathrm{p}}[\kappa T]).$
\subsubsection{Wait-and-see update}
After the decision horizon is completed, the actual load trajectory for decision epoch $\kappa$ is fully realized. This full trajectory $tr^{(N+1)}$ is an additional information and can be used in iterative procedures \eqref{disT} and \eqref{update} to update the augmented tree $\widetilde{\mathcal{T}}^l_{a}[\kappa]$.

The computational results presented in Section \ref{ComRes} bellow, will demonstrate the stability of optimal solutions achieved based on our proposed approach.   
\vspace{-3mm}
\section{Computational Results}\label{ComRes}
We perform our analysis on standard network topology IEEE RTS which is composed of 24 buses, 33 generation units and 38 transmission lines. The data is available in \cite{RTS1979} and also corresponds to IEEE-24 test case in MATPOWER with some add-on scheduling and ramping information. Moreover, in order to test the performance of our approach we will use PJM net-load data from summer 2016 for simulations. The data is mapped across the buses based on the RTS's load profile (cf. \cite{RTS1979}, Table 5).  The algorithms are all implemented in MATLAB. For DMSUC we used YALMIP toolbox and GUROBI 6.5.1 solver. For simulations, all $84$ problem instances (each includes one DMSUC for $\kappa = 0,\ldots, 4$) were parallelized over 32 cores and the numerical experiments were conducted on an Intel Xeon 2.30GHz CPU with 256GB available memory.    
\vspace{-4mm}
\subsection{DMSUC's Optimal Solution}
\begin{figure}
\centering
\includegraphics[width=.98\linewidth]{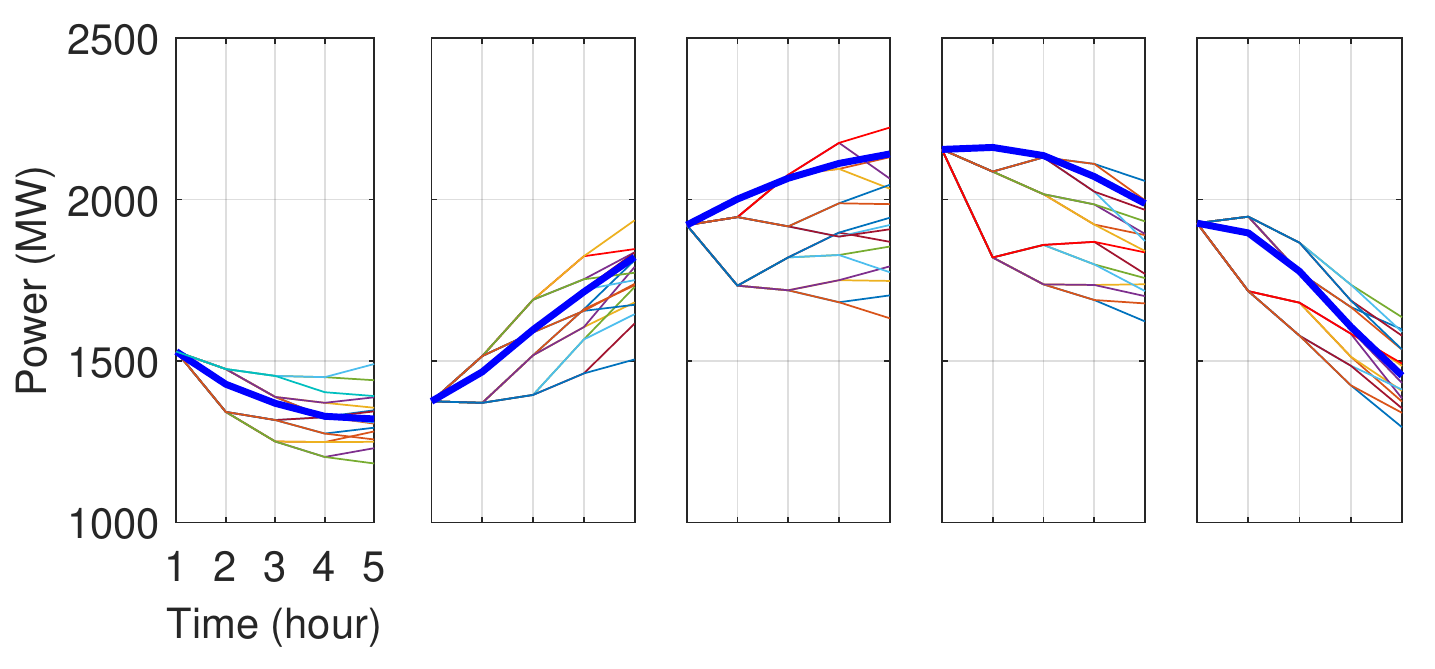}
\caption{\footnotesize{Adopting the correct tree from $\mathfrak{T}^0$ after $\tilde{\xi}_{\mathrm{p}}[\kappa T]$ observed. Actual load (solid blue line) vs. optimal schedule scenarios in the augmented trees.}}\label{ScenVsActualLoad}
\end{figure}

First, we begin by demonstrating how trees $\widetilde{\mathcal{T}}^l[\kappa]$ are chosen from library $\mathfrak{T}^0$ (Fig. \ref{BucketTree}-(b)) and the \textit{hear-and-now} update of the trees 
as is shown in solid blue line in Fig. \ref{ScenVsActualLoad}. Based on \eqref{pickingTrees}, with $\kappa=0,\ldots, 4$, $c=3$, $T = 5$ and the values of actual $\mathrm{p}$resents $\tilde{\xi}_{\mathrm{p}}$, trees 
$\{\widetilde{\mathcal{T}}_a^{1}[0], \widetilde{\mathcal{T}}_a^{2}[1], \widetilde{\mathcal{T}}_a^{3}[2], \widetilde{\mathcal{T}}_a^{3}[3], \widetilde{\mathcal{T}}_a^{2}[4]\}$ are picked from the library and the corresponding root nodes are updated. Consequently,  after solving optimization problem \eqref{cost}-\eqref{state}, Fig. \ref{ScenVsActualLoad}. compares load trajectory with the optimal schedule scenarios for all cohesive decision horizons. From Fig. \ref{ScenVsActualLoad}, it is clear that although the load trajectory is not exactly coinciding with a specific scenario, the optimal schedule scenario tree enfolds the load trajectory almost everywhere and the discrepancy is compensable through ramping actions or procurement of reserves. Moreover, it is probable to observe incidents where actual load  partially lies outside of the range of optimal schedule scenarios as in stage 2 of decision epochs $\kappa = 2,3$ in Fig. \ref{ScenVsActualLoad} and the discrepancy is managed through deployment of reserves.     
Moreover, the average injected power at net-load buses for corresponding consecutive decision horizons is shown in Table. \ref{tab:Average_Pinj}. It is apparent (values in bold face) that with increasing load levels, the direction of injected power at certain buses is subject to change. 
\begin{table}
\centering
\caption{\footnotesize{Ave. injected power at net-load buses for transferring states}}
 \label{tab:Average_Pinj} 
 \begin{tabular}{c | c c c c c} 
 Bus\# &$\widetilde{\mathcal{T}}_a^1[0]$&$\widetilde{\mathcal{T}}_a^2[1]$&$\widetilde{\mathcal{T}}_a^3[2]$&$\widetilde{\mathcal{T}}_a^3[3]$&$\widetilde{\mathcal{T}}_a^2[4]$ \\ [0.75ex] 
\hline
1&-54.4&	-63.4&	\textbf{300.0}&	300.0&	\textbf{-44.8} \\
2&-48.6&	-56.7&	-25.0&	-23.8&	-54.0\\
7&-62.9&	-73.4&	-89.5&	-89.5&	-69.9\\
13&-133.3&	-155.5&	-189.6	&-189.7&	-148.0\\
15&-158.8	&-128.3&	-81.6&	-69.0&	-89.9\\
16&-50.1&	\textbf{47.1}&	82.4&	80.8&	28.1\\
18&232.9&	205.2&	162.4&	162.3&	209.9\\
21&400.0&	400.0&	400.0&	400.0&	380.0\\
22&300.0&	300.0&	300.0&	300.0&	300.0\\
23&330.3&	405.7&	592.1&	580.9&	327.2\\
\hline
\end{tabular}
\end{table}
\vspace{-4mm}
\subsection{Ex-post Performance}
\begin{figure}
\centering
 \includegraphics[width=.98\linewidth,height=4.2cm]{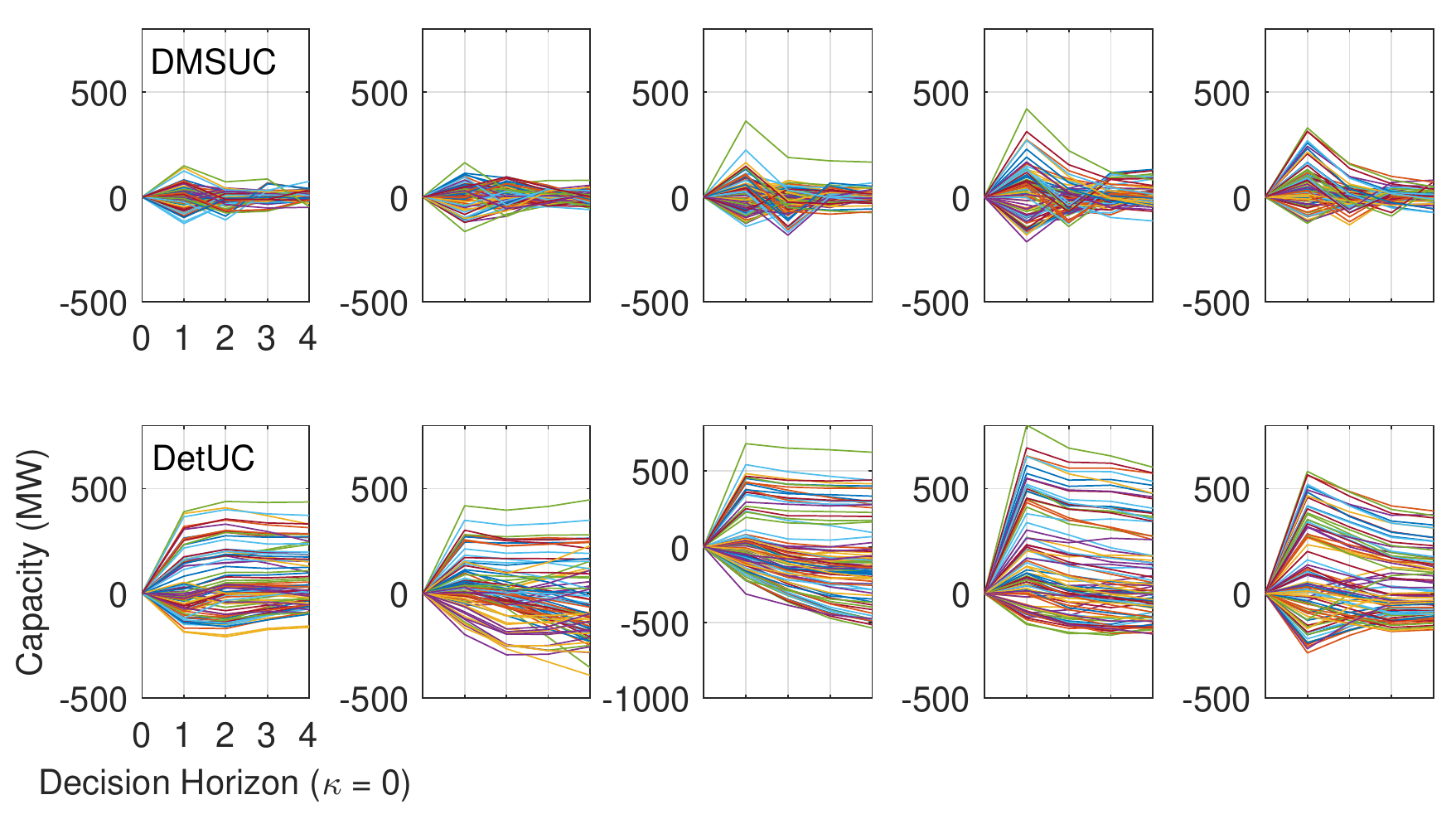}
\caption{\footnotesize{Ex-post schedule difference scenarios for each $\kappa$.}}\label{Ramp_Reserve}
 \end{figure}

\begin{figure}
\centering
 \includegraphics[width=.99\linewidth]{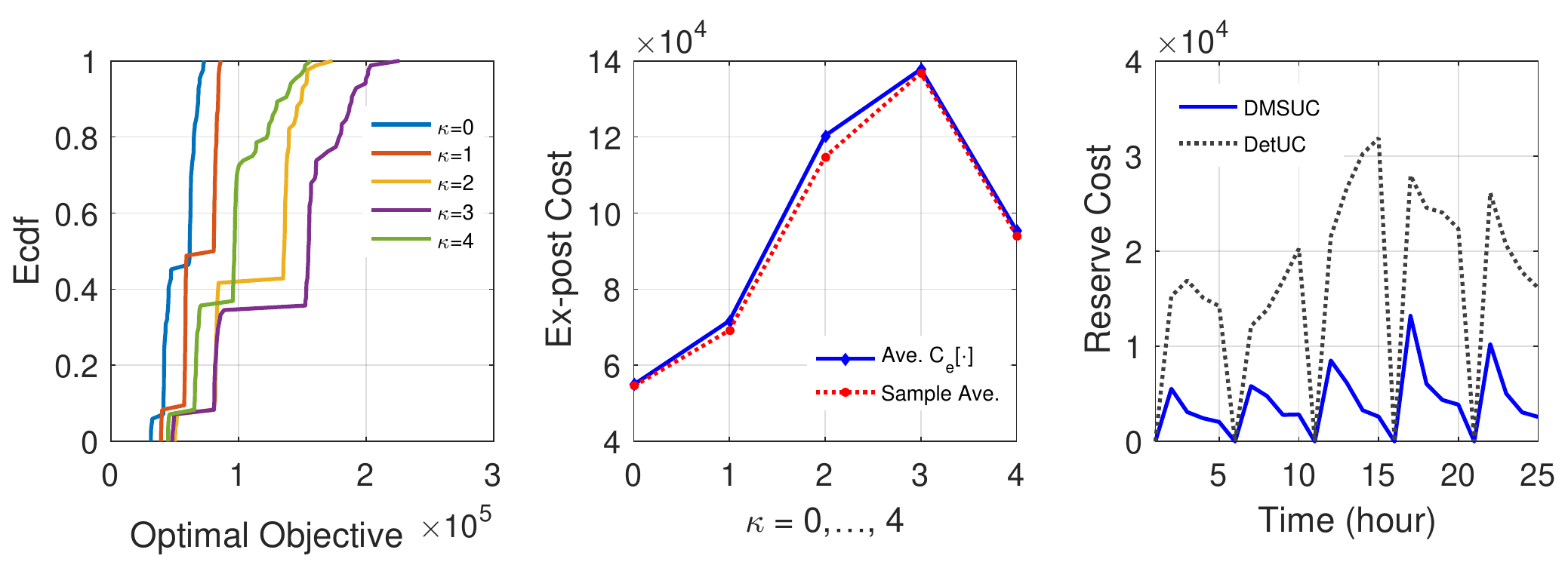}
\caption{\footnotesize{Left: Empirical cdf of DMSUC. Middle: Ave. $C_e[\cdot]$ and Sample Ave. of DMSUC. Right: Reserves cost for residual schedule capacities.}}\label{ExpostCost}
 \end{figure}

In this section we present the performance of DMSUC with respect to schedule scenarios ex-post. We have seen that procedure \eqref{cost}-\eqref{state} is designed to provide optimal scenarios for schedule variables in all decision epochs $\kappa$. By observing the \textit{true} realization of net-load, in or near real time, additional actions might be necessary to adjust the schedule scenario with the real time load. Here, we differentiate between the cost we refer to as \textit{ex-post} and the cost of adjustment actions for procurement of reserves.
\begin{definition}
Ex-post cost $C_{e}[\kappa]$ for each epoch $\kappa$ is defined as:
\vspace{-2mm}
\[C_{e}[\kappa] = \sum_{v\in {\rho^{+}}}\sum_{g\in \mathcal{G}} \left[S^g_{v}+C^g_{v}+\frac{o_{v}^g}{O_n^g}+\frac{d_{v}^g}{O_f^g}\right]\]
with $\rho^{+}$ indicating $\tilde{\xi}_{\rho^{+}}\in \widetilde{\mathcal{T}}^l_a[\kappa]$, which is the closest scenario to the actual load trajectory.
\end{definition}
The provision of reserves' cost, on the other hand, is by finding the difference between the actual load and its nearest optimal schedule path on the tree within that decision horizon. The capacity associated with this distances is an indication for reserve deployment. The corresponding cost is then obtained by taking into account the operational cost of the marginal unit for the required capacity or penalty for over scheduling of the units.  Fig. \ref{Ramp_Reserve}, shows the required adjustment to the schedule scenario for DMSUC and deterministic UC approach. Relying on only a single forecast in deterministic approach and ignoring non-anticipativity of events, results in very conservative schedule scenarios and higher associated costs. Moreover, Fig. \ref{ExpostCost} demonstrates analysis of DMSUC and deterministic UC optimal objective, ex-post and reserves costs for all problem instances corresponding to summer 2016 net-load. On the left, the empirical cdf of the optimal objective is plotted for all epochs $\kappa$ and clearly reveals the daily load pattern. In the middle, the average ex-post cost $C_e[\cdot]$ across all problem instances is plotted and compared to the sample average of DMSUC optimal objective. Alignment of these two measures well promotes the preference of stochastic optimization and performance of our proposed approach, since it does lead to the minimum long term average cost. Finally, on the right, the cost of adjustment actions for both DMSUC and deterministic UC is compared separately for each decision epoch $\kappa$.
\vspace{-2mm}
 \vspace{-1mm}
\section{Remarks and Conclusions}\label{Conclude}
In this paper we presented a new dynamic formulation for MSUC problem and adaptive construction of scenario-tree libraries. In our formulation, the decision horizon of a typical day-ahead UC problem is partitioned into shorter yet consecutive horizons.
This approach not only account for the explicit representation of uncertainty on a scenario tree, but also facilitate the computational complexity of MSUC as mixed-integer MS problem by making it possible to solve a chain of smaller problem instances in a dynamic way. A distinct advantage of our proposed approach comparing to \cite{Guan&etal2016} is that, based on DMSUC approach defining the state variables allows the decision maker to solve MSUS problem in shorter horizons and yet keep the continuity of commitment decisions along the full horizon ($24$ hours) intact. Therefore, redefining the minimum-up and minimum-down times of the generators based on the limitation of decision horizon is eliminated. In addition, when the subset of commitment variables covering the baseline load are fixed at $\kappa = 0$, the size of the optimization problem in further epochs is significantly reduced and only involves updating their states. We also described the algorithm to construct a dictionary of scenario trees by exploiting the cyclostationary property of net-load process. The advantages of constructing the scenario trees from residual processes was also explained. 
Finally, we presented and discussed the performance of DMSUC as opposed to deterministic UC approach. We showed that the ex-post adjustment actions for DMSUC approach results in a profoundly lower reserve costs compared to deterministic UC. In addition, the proximity of ex-post cost and sample average of optimal objective showed yet another strong benefit for utilizing multistage stochastic optimization tool in market environments.        
\vspace{-2mm}
\bibliography{IEEETPS_Paper}
\bibliographystyle{IEEEtran}
\end{document}